\documentclass[11pt]{article}
\usepackage[mathscr]{eucal}

\usepackage{xypic}
\usepackage{amsfonts}
\usepackage{amsmath}
\usepackage{amsthm}
\usepackage{amssymb}
\usepackage{latexsym}
\usepackage{eufrak}

\newtheorem{theorem}[equation]{Theorem}
\newtheorem{proposition}[equation]{Proposition}
\newtheorem{corollary}[equation]{Corollary}

\newtheorem{lemma}[equation]{Lemma}

\newcommand{\beq}{\begin{equation}\label}

\newcommand{\vi}{${\sf {(i)}}\;$}
\newcommand{\vii}{${\sf {(ii)}}\;$}
\newcommand{\viii}{${\sf {(iii)}}\;$}
\newcommand{\viv}{${\sf {(iv)}}\;$}

\newcommand{\ad}{{\mathtt{{ad}}}}

\def\C{{\mathbb{C}}}

\def\R{{\mathsf{R}}}

\def\Z{{\mathbb{Z}}}

\makeatletter
\def\downroundfill{$\m@th \setbox\z@\hbox{$\braceld$}%
  \braceld\leaders\vrule height\ht\z@ depth\z@\hfill\bracerd$}
\def\uproundfill{$\m@th \setbox\z@\hbox{$\braceld$}%
  \bracelu\leaders\vrule height\ht\z@ depth\z@\hfill\braceru$}
\def\overround#1{\mathop{\vbox{\m@th\ialign{##\crcr\noalign{\kern3\p@}
      \downroundfill\crcr\noalign{\kern3\p@\nointerlineskip}
      $\hfil\displaystyle{#1}\hfil$\crcr}}}\limits}
\def\underround#1{\mathop{\vtop{\m@th\ialign{##\crcr
      $\hfil\displaystyle{#1}\hfil$\crcr\noalign{\kern3\p@\nointerlineskip}
      \uproundfill\crcr\noalign{\kern3\p@}}}}\limits}
\makeatother
\begin{document}
\setlength{\parindent}{6mm}
\setlength{\parskip}{3pt plus 5pt minus 0pt}

\centerline{\Large {\textbf{On $m$-quasiinvariants of a Coxeter group}}}

\vskip 4mm
\centerline{\large {\sc {Pavel Etingof and Victor Ginzburg}}}
\vskip 2pt

\begin{abstract}\footnotesize{
Let $W$ be a finite Coxeter group in a Euclidean vector space $V$,
and $m$ a $W$-invariant 
$\Z_+$-valued function on the set of reflections in $W$.
Chalyh and
Veselov introduced in  \cite{CV} an interesting algebra $Q_m$,
called the algebra of {\it $m$-quasiinvariants} for $W$, 
such that $\C[V]^W \subseteq Q_m \subseteq\C[V]$, 
$Q_0=\C[V]$, $Q_m\supseteq Q_{m'}$ if $m\le m'$, and 
$\cap_m Q_m=\C[V]^W$. Namely, $Q_m$ is the
algebra of quantum integrals of the rational Calogero-Moser
system with coupling constants $m$. 
The algebra  $Q_m$ has been studied in \cite{CV}, \cite{VSC}, \cite{FeV}
and  \cite{FV}. In particular, in \cite{FV} Feigin and Veselov  
proposed a number of interesting conjectures 
concerning the structure of $Q_m$, and verified them for dihedral
groups and constant functions $m$. Our goal is to prove some of
these conjectures in the general case.} 
\end{abstract}

\section{Definitions and main results}

We recall some definitions from \cite{FV}. 

Consider a real Euclidean space $V$ of dimension $n$. We will often identify
$V$ and $V^*$ using the inner product on $V$. 

Let $W$ be a finite Coxeter group, i.e. a finite group 
generated by reflections 
of $V$. Let $N=|W|$. Let $\Sigma$ 
be the set of reflections in $W$, and $\Pi_s$ be the reflection
hyperplane for a reflection $s$. Let $m:\Sigma\to \Z_{\ge 0}$, $s\mapsto m_s$, be a
$W$-invariant function (called the multiplicity function). 
A complex polynomial $q$ on $V$ is said to be an $m$-quasiinvariant 
(under $W$) if, for each $s\in\Sigma$, the function $x\mapsto 
q(x)-q(sx)$ vanishes up to order $2m_s+1$ at the
hyperplane $\Pi_s$. Such polynomials form a graded subalgebra in
the graded algebra
$\C[V]=\bigoplus_{i\ge 0}\;\C[V]\langle i\rangle$, which will be denoted by $Q_m$. It is obvious 
that $Q_m$ contains as a subalgebra the ring $\C[V]^W$ of
invariant polynomials. We denote by $I_m$ the ideal in $Q_m$
generated by the augmentation ideal in $\C[V]^W$. 
This is a graded ideal in $Q_m$. 

The following two theorems, conjectured in \cite{FV}, 
are two of the main results of this paper. 

Let $T$ be any graded
complement of $I_m$ in $Q_m$. 

\begin{theorem}\label{free} $Q_m$ is a free
  module over $\C[V]^W$, of rank $N$. More specifically, the
  multiplication mapping defines a graded isomorphism 
$\C[V]^W\otimes T\to Q_m$. In particular,
  $\text{dim}(Q_m/I_m)=\text{codim}(I_m)=N$.  
\end{theorem}

Consider now the $N$-dimensional graded algebra $R_m=Q_m/I_m$. 
Let $d=\sum_{s\in \Sigma}{(2m_s+1)}$. 

\begin{theorem}\label{gor}
\vi The space $R_m\langle d\rangle$ is
  one dimensional. 

\vii {\sf (Poincare duality).}\, 
The multiplication mapping 
$R_m\langle j\rangle\times R_m\langle d-j\rangle\to R_m\langle d\rangle$ 
is a nondegenerate pairing 
for any $j$. In particular, 
the Poincare polynomial 
$P_{R_m}(t)$ is a palindromic polynomial
of degree $d$ (i.e. $P_{R_m}(t^{-1})=t^{-d}P_{R_m}(t)$), and the algebra 
$R_m$ is Gorenstein. 

\viii The algebra $Q_m$ is Gorenstein.
\end{theorem}

The proofs of Theorem \ref{free} and Theorem \ref{gor}
are given in the next few sections. 

{\bf Remarks.}  \vi For $m=0$, the quasiinvariance condition is
vacuous, so $Q_m=\C[V]$. 
Thus, for $m=0$ Theorem \ref{free} 
reduces to the Chevalley theorem, which claims 
that $\C[V]$ is free over $\C[V]^W$. 
Therefore, Theorem \ref{free} is an m-version of the Chevalley
theorem. (We note, however, that our proof of Theorem \ref{free}
makes use of the Chevalley theorem, so we do not obtain a new proof
of the Chevalley theorem). Theorem \ref{gor} for $m=0$ is also
well known, and is due to Steinberg. 

\vii If $W$ is a Weyl group, this
theorem has a topological interpretation, since in this case 
$R_m$ is the cohomology algebra of the flag variaty for the
correponding complex semisimple Lie group. 

\medskip

{\bf Acknowledgments.}{ \footnotesize{The work of the first author was partly
conducted for the Clay Mathematics Institute, and was partially
supported by the NSF grant DMS-9988796.   
We are grateful to E.Opdam, who helped us with the proof of
Theorem \ref{00}. We are indepted 
to M.Feigin,  G.Felder, and A.Veselov for useful
discussions and for making their results
available before publication. In
particular, we note that at the time this paper was written,
the preprint \cite{FeV} had not yet appeared, and that we owe all 
our information about the results of \cite{FeV} to private
communications with the authors. 
Finally, the first author is grateful to J.Starr for
explanations about Gorenstein algebras.}}

\section{The Calogero-Moser quantum integrals}

Let us recall some known facts about the quantum Calogero-Moser
systems. 

Let $\mu$ be a complex valued 
$W$-invariant function on $\Sigma$.
Let $H=H(\mu)$ be the differential operator
$$
H(\mu)=\Delta_V-\sum_{s\in
  \Sigma}\frac{2\mu_s}{(\alpha_s,x)}\partial_{\alpha_s},
$$
where $\alpha_s\in V^*$ is an eigenvector for $s$ with eigenvalue
$-1$, and ${\partial}_\alpha$, $\alpha\in V^*$, is the derivation of $\C[V]$ 
defined by ${\partial}_\alpha(\beta)=(\alpha,\beta)$, $\beta\in V^*$. 
The operator $H(\mu)$ is homogeneous, and has degree $-2$. 
It is called the Calogero-Moser quantum hamiltonian. 

The following result shows that for any $\mu$, the Calogero-Moser 
hamiltonian defines a quantum integrable system. 

Let $p_i$, $i=1,...,n$, be a set of 
algebraically independent homogeneous generators of
$\C[V]^W$, and $d_i$ be their degrees. 

\begin{theorem} \cite{OP},\cite{Ch} For any $i$, there exists 
a differential operator $L_i=L_i(\mu)$ of homogeneity
degree $-d_i$, with principal
symbol $p_i(\xi)$, such that $[L_i,H]=0$. 
The operators $L_i$ are 
regular outside of the 
reflection hyperplanes and algebraically independent. \qed
\end{theorem}

In fact, the operator $L_i$ with the above properties is unique.
It is obtained by evaluating the polynomial $p_i$ on the Dunkl
operators, and then restricting the resulting operator to
invariant functions. 
The algebra generated by the operators $L_i$ will be denoted by $A_\mu$.   

Let $D$ be the algebra of differential operators on $V$ with
rational coefficients, which are regular outside of the
reflection hyperplanes. 
For generic $\mu$, the algebra $A_\mu$ is a maximal commutative
subalgebra of $D$. However, for integer $\mu$, additional quantum
integrals turn out to exist. They are described by the theorem below
which follows from \cite{CV},\cite{VSC};
see also \cite{FV}.

\begin{theorem} Let $m$ be a nonnegative
  integer invariant function on $\Sigma$.
There exists a unique algebra homomorphism 
$\phi: Q_m\to D$ such that $\phi(p_i)=L_i(m)$, and the symbol 
of the operator $\phi(q)$ is $q$ for any homogeneous polynomial
$q$. This homomorphism maps elements of degree $d$ to elements of
degree $-d$. The image of $\phi$ is the centralizer 
of $A_m$.  \qed
\end{theorem}

We will denote $\phi(q)$ by $L_q(m)$ or shortly by $L_q$. 
In particular, $L_i=L_{p_i}$.  
Thus, if $q$ is $W$-invariant, then 
$L_q(\mu)$ makes sense not only for positive integer, but actually
for any complex $\mu$. 

{\bf Remark.} It is worth mentioning an explicit formula for
$L_q$, due to Berest: if $q$ is homogeneous of degree $r$ then 
$L_q=c(\ad H(m))^rq$, where $c$ is a constant. 

Now let us introduce the $\psi$-function, which plays the main
role in this paper. 

\begin{theorem}(\cite{CV},\cite{VSC}) 
There exists a unique, up to scaling, function 
$\psi_m(k,x)$ on $V\times V$, of the form 
$P(k,x)e^{(k,x)}$, where $P$ is a polynomial, such that 
for any $q\in Q_m$ one has $L_q^{(x)}(m)\psi_m(k,x)=q(k)\psi_m(k,x)$. \qed
\end{theorem}

The function $\psi_m$ is called the Baker-Akhiezer function.
We will denote it simply by $\psi$, assuming that 
$m$ is fixed.  
It is clear that $\psi$ is homogeneous in the sense
$\psi(tk,x)=\psi(k,tx)$. One can also show that the highest term 
of $P(k,x)$ is proportional to $\delta_m(k)\delta_m(x)$, where 
$\delta_m(x)=\prod_{s\in \Sigma}(\alpha_s,x)^{m_s}$
(up to scaling, it is independent on the choice of $\alpha_s$). 

\begin{theorem}(\cite{CV},\cite{VSC}) The function $\psi$ is symmetric 
under the interchange of $k$ and $x$. In particular, it has the
bispectrality property: $L_q^{(k)}\psi(k,x)=q(x)\psi(k,x)$.  \qed
\end{theorem}

\begin{theorem} (\cite{CV},\cite{VSC}) The function $\psi(k,x)$ is
  $m$-quasiinvariant with respect to both variables.  \qed
\end{theorem}

The following result plays a key role in this paper. 

\begin{theorem}\label{00} $\psi(0,0)\ne 0$. 
\end{theorem}

\begin{proof} According to \cite{DJO}, integer valued
  multiplicity functions are nonsingular in the sense of 
\cite{DJO}. This means 
(see \cite{DJO}, p. 247) that there exists a generalized Bessel
function $J_m(k,x)$, a holomorphic $W$-invariant in both variables
solution of the system of differential equations
$$
L_i^{(x)}(m)J_m(k,x)=p_i(k)J_m(k,x)
$$
 This function is unique up to 
scaling, and can be normalized by the condition $J_m(0,0)=1$. 

Now consider the function $K_m(k,x)=\sum_{w\in W}\psi(k,wx)$. 
This function is a holomorphic invariant solution of the above
system, so it must be proportional to $J_m$. 
This implies that $\psi(0,0)\ne 0$, as desired. 
\end{proof} 

In fact, there is an explicit product formula for
$\psi(0,0)$, for the normalization of the highest term of 
$P$(= the polynomial factor in the the $\psi$-function)
 to be $\delta_m(x)\delta_m(k)$, with all 
$\alpha_s$ having squared length $2$. Such a formula can be
deduced from \cite{DJO}. However, we will not discuss 
this formula, and will choose the normalization of $\psi$ 
such that $\psi(0,0)=1$.  

\section{The pairing on $Q_m$}

Let us expand $\psi(k,x)$ into a Taylor series. 
Since $\psi$ is $m$-quasiinvariant with respect to both variables, 
 we can consider $\psi(k,x)$ as an element of $Q_m\hat\otimes
 Q_m$, where $\hat\otimes$ is the completed tensor product. 
Furthermore, because of the homogeneity of $\psi$, we have
$\psi=\sum_{j\ge 0}\psi^{(j)}$, where 
$\psi^{(j)}\in Q_m\langle j\rangle\otimes Q_m\langle j\rangle$. 

\begin{proposition} $\psi$ is a nondegenerate tensor, i.e. its
left (or right) tensorands span $Q_m$. 
In other words, the tensorands of $\psi^{(j)}$ span 
$Q_m\langle j\rangle$ for all $j\ge 0$. 
\end{proposition}

\begin{proof} Let $Q_m'\subset Q_m$ be the span of left (or right) 
tensorands of $\psi$. This is a graded subspace of $Q_m$. 
Let $q_i$ be a homogeneous basis of $Q_m'$. 
Then we can write $\psi$ in the form $\sum_i q_i(k)q^i(x)$, where
$q^i(x)$ is another homogeneous basis of $Q_m'$. Thus, for any 
$q\in Q_m$, we have 
$$
q(k)\psi(k,x)=(L_q^{(x)}\psi)(k,x)=\sum_i q_i(k)(L_qq^i)(x).
$$
But the function $q(k)\psi(k,x)$ is
analytic. Thus, $L_qq_i$ cannot have poles and hence is a polynomial. 

Let us now substitute $x=0$ in the last equality. 
Since $\psi(k,0)=\psi(0,0)=1$, we get 
$q(k)=\sum a_iq_i(k)$, where $a_i=(L_qq_i)(0)$. 
This sum is clearly finite. Thus, $q\in Q_m'$, i.e. 
$Q_m'=Q_m$, as desired. 
\end{proof}

This proposition 
and the fact that $\psi$ is an eigenfunction of $L_q$ 
has the following corollary, which 
is also proved in \cite{FV} by another method:

\begin{corollary}\cite{FV} 
For any $q\in Q_m$ one has $L_q(Q_m)\subset
  Q_m$. 
\end{corollary}

Consider now the symmetric bilinear form on $Q_m$ inverse to 
the element $\psi$. This form is nondegenerate.
We will denote it by $(p,q)_m$, or simply by $(p,q)$. 
The next theorem summarizes the properties of this form. 

\begin{theorem}

\vi  $(,)$ is $W$-invariant, and 
$(p,q)=0$ if $p,q$ are homogeneous of different degrees. 

\vii  $(p,q)=(L_qp)(0)$. 

\viii  $(p,q)=L_p^{(x)}L_q^{(k)}\psi(x,k)|_{x=k=0}$

\viv $L_q^*=q$ under $(,)$. 

\end{theorem}
 
\begin{proof} (i) is clear. 

Proof of (ii): Let $q_i$ be a homogeneous basis of $Q_m$, 
and $q^i$ the dual basis. Then $\psi(k,x)=\sum_iq_i(k)q^i(x)$. 
Applying $L_{q_l}^{(x)}$ to this equation, we get 
$$
q_l(k)\psi(k,x)=\sum_i q_i(k)(L_{q_l}q^i)(x)
$$
Substituting $x=0$, we get 
$$
q_l(k)=\sum_i q_i(k)(L_{q_l}q^i)(0),
$$
so $(L_{q_l}q^i)(0)=\delta_l^i$, as desired. 

Proof of (iii):  By (ii), the right hand side is 
$\sum_i(p,q^i)(q,q_i)=(p,q)$. 

Proof of (iv): by (ii) $(L_qp_1,p_2)=(L_{p_2}L_qp_1)(0)=
(L_{qp_2}p_1)(0)=(qp_2,p_1)$. 
\end{proof}

The above results on the form on $Q_m$ imply the following.
Let $D_m$ be the algebra generated by $q,L_q$, $q\in Q_m$. 

\begin{proposition} $Q_m$ is an irreducible $D_m$-module. 
\end{proposition}

\begin{proof}
First of all, $D_m$ clearly contains the Euler vector field, so any 
submodule of $Q_m$ has to be graded. Thus, it is sufficient 
to show that for any homogeneous element $q\in Q_m$, one has
$q\in D_m1$ and $1\in D_mq$. 
But this is clear, 
since for any homogeneous $q\in Q_m$ one has $q=q1$,
  and $1=L_pq$ for $p$ of degree $\text{deg}(q)$ 
such that $(p,q)=1$ (which exists by
nondegeneracy of the form).
\end{proof} 

\section{Proof of Theorem \ref{free}} 

We are ready to prove Theorem \ref{free}. For this purpose, for
any $k\in V$, define the subspace $H_m(k)$ of the power series 
completion $\C[[V]]$ of $\C[V]$, which consists of solutions 
of the differential equations $L_if=p_i(k)f$, 
$i=1,...,n$. 

\begin{theorem} (essentially, contained in \cite{FV}) 

\vi  $\text{dim}H_m(k)=N$ for all
  $k$. 

\vii  $H_m(k)$ is contained in the power series completion 
$\hat Q_m$ of $Q_m$. 
 
\viii  $H_m(0)$ is a graded subspace of $Q_m$.
\end{theorem}

\begin{proof}  
(i) Looking at the symbols of $L_i$ and using
  the Chevalley theorem, we conclude that the dimension cannot be 
more than $N$ (since this is the dimension of the space of ``abstract'' 
solutions of the system, in the sense of differential Galois theory).
On the other hand, for generic $k$, it is easy to see 
that the functions $\psi(x,wk)$, $w\in W$, are linearly
  independent elements of $H_m(k)$. Thus the dimension is
  generically (and hence always) greater than or equal to $N$. 
Combining the two results, we get that the dimension 
is exactly equal to $N$. 

(ii) The statement says that elements of $H_m(k)$ are
$m$-quasiinvariant. This is clear for generic $k$ since 
we showed in the proof of (i) that $\psi(x,wk)$ is a 
basis of $H_m(k)$. Therefore, it is true for all $k$. 

(iii) This is clear, as $H_m(0)$ is graded and finite
dimensional. 
\end{proof}

\noindent
{\bf Remark.} 
The space $H_m=H_m(0)$ is called in \cite{FV} the space of 
$m$-{\it harmonic} polynomials. 

Now let $I_m(k)$ be the ideal in $Q_m$ generated by 
the polynomials $p_i-p_i(k)$. In particular, $I_m(0)=I_m$. 

\begin{lemma} $I_m(k)$ is the orthogonal complement 
of $H_m(k)$ in $Q_m$ with respect to the form $(,)$.
\end{lemma} 

\begin{proof}
This follows from the fact that $L_q^*=q$. 
\end{proof}

\begin{corollary} $Q_m/I_m(k)=H_m(k)^*$ is a flat family of
  vector spaces, of dimension $N$. 
\end{corollary}

Now let us consider $Q_m$ as a module over $\C[V]^W$. The fiber 
of this module at the point $k$ is $Q_m/I_m(k)$. Since this family 
is flat, the module $Q_m$ is locally free. 
But since it is graded, it is freely generated by any 
local homogeneous generators 
$t_1,...,t_{N}$ at the point $k=0$. This proves Theorem \ref{free}.

The proof of Theorem \ref{free} we gave also implies the
following Corollary, which was conjectured in \cite{FV}:
 
\begin{corollary} \label{poin} One has the following identity for 
Poincare series: 
$$
P_{Q_m}(t)=\frac{P_{H_m}(t)}{\prod_{i=1}^n(1-t^{d_i})}.
$$
In particular, $P_{R_m}(t)=P_{H_m(t)}$. 
\end{corollary} 

The polynomial $P_{H_m}$ is
calculated
in \cite{FeV}. Thus, the Corollary allows one to compute
$P_{Q_m}$
and $P_{R_m}$. 

\section{Some determinants}

In this section we will calculate the order of vanishing of some 
determinants, which will be used later. 

Let $\delta_{2m+1}(x)=\prod_s \alpha_s(x)^{2m_s+1}$ be the
m-version of the discriminant. 

\begin{proposition}\label{wedge} Let $k\in V$.
The polyvector $u(k)=\bigwedge_{g\in W}\psi(gk,*)\in \bigwedge^N
H_m(k)$ is of the form $\delta_{2m+1}(k)u_*(k)$, where $u_*(k)$
is a nonvanishing section of the line bundle $\Lambda^N H_m(k)$
over $V$. 
\end{proposition}

\begin{proof} 
It is clear that when $k$ is regular then $u(k)\ne 0$, because in
this case $\psi(gk,x)$ form a basis of $H_m(k)$. Thus, it is
sufficient to show that 
$u(k)$ has exactly the prescribed order of vanishing on the
reflection hyperplanes. 

To show this, 
let $k_0$ be a generic point of $\Pi_s$, and 
$v$ be a nonzero vector orthogonal to $\Pi_s$. 
Define the function 
$$\tilde \psi(k_0,x)=\lim_{\epsilon\to 0}\frac{\psi(k_0+\epsilon
  v,x)-\psi(k_0-\epsilon v,x)}{\epsilon^{2m_s+1}}.$$
This function is obviously nonzero, 
and is well defined up to normalization (the
normalization depends on $v$). 

It is easy to see that the functions 
$\psi(gk_0,x)$ and $\tilde \psi(gk_0,x)$, $g\in W/(1,s)$,
are linearly independent, and form a basis of $H_m(k_0)$. 
Therefore, the wedge product 
$$
u'(k)=\bigwedge_{g\in W/(1,s)}(\psi(gk,x)+\psi(gsk,x))
\wedge
\bigwedge_{g\in W/(1,s)}\frac{\psi(gk,x)-\psi(gsk,x)}{\alpha_s(k)^{2m_s+1}}
$$ 
has a nonzero finite limit as $k\to k_0$. But it is clear that $u'(k)$
is a constant multiple of $u(k)/\alpha_s(k)^{N(2m_s+1)/2}$. This
implies the required statement.
\end{proof}

Let $T$ be a graded linear complement to $I_m$ in $Q_m$. 
Let $t_i$, $i=1,...,N$, be a 
homogeneous basis of $T$. Let $k\in V$. 
Let $A(k)$ be the matrix whose entries are $t_i(gk)$,
$i=1,...,N$, $g\in W$. 

\begin{lemma}\label{det1}
$$
\det A(k)=c\delta_{2m+1}(k)^{N/2},
$$
where $c$ is a nonzero constant. 
\end{lemma}

\begin{proof}
First, note that the pairing $(,): T\times H_m(k)\to \C$ is
nondegenerate for any $k\in V_\C$. Indeed, since $T$ is graded, the degeneracy
locus of this pairing in the k-space is invariant under dilations.
Also, this locus is clearly closed. So, if it is nonempty, 
it must contain zero. 
But the pairing between $T$ and $H_m(0)$ is nondegenerate by the
definition of $T$. 

This implies that for any regular point $k\in V_\C$, 
the evaluation map $T\to \C[W\cdot k]$ is 
an isomorphism (since for $f\in T, f(gk)=(f,\psi(gk,*))$, 
and $\psi(gk,*)$ is a basis of $H_m(k)$). 
Thus, $\det A(k)$ is nonzero outside of 
the reflection hyperplanes in $V_\C$. So it suffices to check 
that $\det A(k)$ has exactly the predicted degree of vanishing on the
hyperplanes, i.e. degree $N(2m_s+1)/2$ on $\Pi_s$.  

Let us first check that 
$\det A(k)$ has degree of vanishing at least $N(2m_s+1)/2$ on
$\Pi_s$. To this end, look at the limit in which $k$
approaches a generic point $k_0$ on a hyperplane $\Pi_s$. 
Since $t_i$ are quasiinvariants, for any $g\in W/(1,s)$, the difference 
between $t_i(gk)$ and $t_i(gsk)$ is of the order at least 
$\alpha_s(k)^{2m_s+1}$ in this limit. This gives the desired lower
bound. 

Now let us obtain the upper bound. 
As we mentioned, the pairing \linebreak $T\times H_m(k)\to
\C$ given by $(,)_m$ is nondegenerate. 
Thus, there exists a basis $f_j=f_j^{(k)}$ of 
$H_m(k)$ such that $(L_{t_r}f_j)(0)=\delta_r^j$. 

Let us express the solutions $\psi(gk,x)$ via this basis. 
It is clear that $\psi(gk,x)=\sum_j c_j(g)f_j(x)$, 
where $c_j(g)=L_{t_j}^{(x)}\psi(gk,x)|_{x=0}=
t_j(gk)$. Thus, 
$$
\psi(gk,x)=\sum_j t_j(gk)f_j(x),
$$
and
$$
u(k)=\det A(k)\cdot \bigwedge_{j=1}^Nf_j^{(k)}. 
$$
The second factor is holomorphic in $k$. 
Thus, the lower upper bound follows from Proposition \ref{wedge}.
The Lemma is proved. 
\end{proof} 

\begin{corollary} (proved also in \cite{FeV})
The number $P_{H_m}'(0)=\sum \text{deg}(t_j)$ equals to \linebreak
$(N/2)\sum_s(2m_s+1)$.
\end{corollary}

\begin{proof}
This is obtained from Lemma \ref{det1}
by comparing the degrees of the two sides.
\end{proof} 

\section{Linear independence theorem and [FV]-conjectures}

Let $T,t_i$ be as in the previous section. 
Recall \cite{FV} that $\delta_{2m+1}(x)\in H_m$. 
Hence for any $j$, one has $L_{t_j}\delta_{2m+1}\in H_m$. 

\begin{theorem}\label{linindep} 
The elements $L_{t_i}\delta_{2m+1}$ are linearly
  independent, 
and hence form a basis of $H_m$. 
\end{theorem}

\begin{proof} Let $k\in V$ be regular. 
Consider the function 
$\delta_{2m+1}^{(k)}(y)=\frac{\sum_{h\in
    W}(-1)^h\psi(hk,y)}
{\delta_{2m+1}(k)}$. It is easy to see from quasiinvariance of
$\psi$ that 
this function (as a function of $k$)
extends to a holomorphic function on $V_\C$. 
In particular, there exists a limit 
$$
\delta_{2m+1}^{(0)}(y):=
\lim_{k\to 0}\delta_{2m+1}^{(k)}(y),
$$
and $\delta^{(0)}_{2m+1}(y)$ is an antisymmetric 
m-harmonic polynomial. Hence, $\delta^{(0)}_{2m+1}(y)=b\delta_{2m+1}(y)$, 
where $b\in\C$.

Consider the polyvector
$$
\Delta(k)=\bigwedge_{j=1}^NL_{t_j}\delta_{2m+1}^{(k)}=
\bigwedge_{j=1}^NL_{t_j}\frac{\sum_{h\in
    W}(-1)^h\psi(hk,*)}{\delta_{2m+1}(k)}
$$
(the last expression applies to regular $k$ only).
We have $\Delta(0)=b^N\bigwedge_{j=1}^NL_{t_j}\delta_{2m+1}$. Thus, 
it is sufficient for us to show that $\Delta(0)$ is nonzero. 

For regular $k$, we have  
$$
\Delta(k)=\bigwedge_{j=1}^N\frac{\sum_{h\in
    W}(-1)^ht_j(hk)\psi(hk,*)}{\delta_{2m+1}(k)}=
\delta_{2m+1}(k)^{-N}\cdot \det (A(k)J)\cdot u(k),
$$
where $J_{hh'}=(-1)^h\delta_{hh'}$. 
Thus, by the lemmas on determinants, 
$$
\Delta(k)=\pm \delta_{2m+1}(k)^{-N}\delta_{2m+1}(k)^{N/2}
(\delta_{2m+1}(k)^{N/2}u_*(k))=\pm u_*(k).
$$
But we have seen above that $u_*(0)\ne 0$. 
Hence, $\Delta(0)\ne 0$,
as desired.
\end{proof}

\noindent
{\bf Remark.} In particular, we have shown that $b\ne 0$. 

\begin{corollary} \label{con1} (Conjecture 1 of \cite{FV}). 
Consider the linear map $\pi_m: Q_m\to H_m$, given by 
$\pi_m(q)=L_q\delta_{2m+1}$.  
Then $\pi_m$ is surjective, and 
the kernel of $\pi_m$ is $I_m$.
\end{corollary}

\begin{proof} The first statement is clear
from Theorem \ref{linindep}. 
The second statement follows from the first one, since
$\text{Ker}(\pi_m)\supset I_m$, 
and $\text{codim}(I_m)=\text{dim}(H_m)$. 
\end{proof}

Consider now the bilinear form $<,>:Q_m\times Q_m\to \C$ defined
by $<p,q>=(p,\pi_mq)=(L_{pq}\delta_{2m+1})(0)$. 
It is clear that this form is symmetric, and its kernel 
contains $I_m$. Thus, this form induces a form 
$<,>$ on the algebra $R_m=Q_m/I_m$, which has homogeneity degree 
$d=\sum_s(2m_s+1)$. 

\begin{proposition} The form $<,>$ on 
$R_m$ is nondegenerate. 
\end{proposition}

\begin{proof} It is sufficient to show that the restriction of
  $<,>$ to $T$ is nondegenerate. But this follows from the fact
  that $\pi_m:T\to H_m$ is an isomorphism (Corollary \ref{con1}), and that 
$(,):T\times H_m\to \C$ is nondegenerate (definition of $T$).  
\end{proof}

\vskip .1in
{\bf Proof of Theorem \ref{gor}.}
(i) Since the form $<,>$ has degree $d$ and is nondegenerate, we
have $\dim R_m\langle j\rangle=\dim R_m\langle d-j\rangle$. In particular, 
$\dim R_m\langle d\rangle=1$. 

(ii) It is clear that 
$R_m\langle d\rangle$ is spanned by the image of $\delta_{2m+1}$. Indeed this 
image is clearly nonzero (as $\delta_{2m+1}$ is the lowest degree
antisymmetric element in $Q_m$, see \cite{FV}), which implies 
that it spans $R_m\langle d\rangle$. 

Now, it is easy to see that the multiplication mapping in question 
 $p,q\mapsto p*q$ is 
proportional to $p,q\mapsto <p,q>\delta_{2m+1}$. The 
nondegeneracy conclusion follows, and the Gorenstein property
 follows from nondegeneracy. 

(iii) The algebra $Q_m$ is graded and is free as a module over 
$\C[p_1,...,p_n]$. By standard results of commutative algebra
(see \cite{Eis}, Chapter 21), 
this implies that $Q_m$ is Gorenstein if and only if so 
is $Q_m/(p_1,...,p_n)$. But $(p_1,...,p_n)=I_m$, so 
$Q_m/(p_1,...,p_n)=R_m$, and we know from (ii) that $R_m$ is
Gorenstein. The theorem is proved. 
\qed

\vskip .1in

We conclude this section with an 

\noindent
{\bf Alternative proof of Theorem \ref{gor}}.\quad
This proof is based on the  following remarkable result
due to R. Stanley \cite{St}

\begin{theorem}\label{stan}
A positively graded Cohen-Macaulay domain is Gorenstein 
if and only if its Poincare series $h(t)$ is a rational function 
which satisifies the equation 
$h(t)=(-1)^nt^lh(t^{-1})$ for some $l$ and for $n$ being the 
(algebro-geometric) dimension 
of $A$. \qed
\end{theorem}

Let us use this result to prove Theorem \ref{gor}. 
First of all, we note that by Theorem \ref{free}, the algebra 
$Q_m$ is Cohen-Macaulay (since it is a free module over a smooth
subalgebra, see \cite{Eis}, Corollary 18.17). It is also
positively graded
and does not have zero divisors (as it is a subring of $\C[V]$). 

Next, we cite a result of \cite{FeV}:

\begin{proposition} The polynomial $P_{H_m}$ is 
a palindromic polynomial of degree $d$. That is, 
$P_{H_m}(t^{-1})=t^{-d}P_{H_m}(t)$. 
\end{proposition}

Thus, the same is true about $P_{R_m}$, since by 
Corollary \ref{poin}, $P_{R_m}=P_{H_m}$. 
Therefore, it is easy to check that Stanley's criterion is
satisfied, and by Theorem \ref{stan}, 
$Q_m$ and hence $R_m$ are Gorenstein
algebras. This proves Theorem \ref{gor}. 

The other results of this section 
follow easily from this. Indeed, since $p*q=<p,q>\delta_{2m+1}$, 
$p,q\in R_m$, we get that $<,>$ is nondegenerate, and since 
$<p,q>=(p,\pi_mq)$, we get that $\pi_m: R_m\to H_m$ is an
isomorphism.  

\medskip

{\bf Remark.} One can define an obvious 
analog $Q_m(\Sigma)$ of the algebra $Q_m$ 
for any arrangement of hyperplanes $\Sigma$ in a Euclidean space,
and any positive
integer function $m$ on $\Sigma$. However, in general this algebra
will not be as nice as $Q_m$. 

For example, suppose that $\Sigma$ consists 
of two lines through $0$ in the plane, and $m=1$. 
If the lines are perpendicular, we have the Coxeter configuration
for the group $W=(\Z/2)^2$, so $Q_m$ has Poincare 
series $P(t)=(\frac{1-t+t^2}{1-t})^2$ and is Gorenstein. 
However, if the lines are not perpendicular, it is not difficult
to show that 
the Poincare series of $Q_m(\Sigma)$ is given by $\hat
P(t)=P(t)-t^2$, so 
$$
\hat P(t)=\frac{1-2t+2t^2}{(1-t)^2}.
$$
It is clear that this function does not satisfy Stanley's
criterion. Therefore, 
we see that $Q_m(\Sigma)$ is not Gorenstein unless the lines are 
perpendicular, i.e. unless we have a Coxeter configuration. 

It
would be interesting to know whether this phenomenon occurs for more
general classes of configurations. 
\bigskip

\section{A counterexample}

We have only proved a part of conjectures from \cite{FV}. 
The rest of the conjectures claim that

\vi  (Conjecture 2) the restriction of the map $\pi_m$ to 
$H_m$ is an isomorphism $H_m\to H_m$. 

\vii  (Conjecture 3) $Q_m$ is generated by $H_m$ over $\C[V]^W$. 

\viii  (Conjecture 2*) The restriction of the pairing $<,>$ to
$H_m$ is nondegenerate. 

These conjectures were proved in \cite{FV} for 
dihedral groups and constant functions $m$. 

Unfortunately, it turned out that 
these conjectures do not hold for general $W$ and $m$. 
This is demonstrated by the following example. 

Let $W$ be of type $B_6=C_6$, so the roots are $\pm e_i$, $\pm
e_i\pm e_j$, and the basic invariant polynomials 
are $p_j=\sum_{i=1}^6x_i^{2j}$, $j=1,...,6$.
 Let $m=1$ for the short roots $\pm e_i$, and 
$m=0$ for the long roots $\pm e_i\pm e_j$. 

We claim that (i),(ii), and (iii) are not satisfied for these $W$
and $m$. To see this, we 
let $M$ be the operator $\frac{d^2}{dx^2}-(2/x)\frac{d}{dx}$. 
Then the operators $L_{p_i}$ are 
$$
L_{p_1}=M_1+...+M_6,...,
L_{p_6}=M_1^6+...+M_6^6
$$
(where $M_i$ is $M$ acting in the variable $x_i$) 

It is clear that the polynomial $u=x_1^3$ is m-harmonic. 
It is easy to check that the polynomial
$$
v=x_1^3(x_1^2+...+x_6^2)
$$ 
is also m-harmonic. 

But $u,v$ are linearly dependent over $\C[V]^W$. 
Thus, since $Q_m$ is free over $\C[V]^W$ of rank $N$, and 
$H_m$ has dimension $N$, it is impossible that $H_m$ generates
$Q_m$ over $\C[V]^W$, which disproves (ii). 

We also have $H_m\cap I_m\ne 0$. Indeed, consider a lowest
degree element $q$ in the orthogonal complement to $\C[V]^WH_m$ in
$Q_m$ with respect to $(,)$. The elements $L_{p_i}q$ must also be in this 
complement, so $L_{p_i}q=0$ and $q\in H_m$. On the other hand,
$(H_m,q)=0$, so 
$q\in I_m$. This disproves (i). 

Finally, since $<p,q>=(p,\pi_m q)$ on $H_m$, (iii) fails since
$\pi_m$ must have nonzero kernel in $H_m$.  

\section{The shift operator}

It is interesting to point out the relation of the above with
the shift operator. 

The following theorem is due to Opdam (see \cite{Op}). 

Let $\mu: \Sigma\to \C$ be an invariant function. 

\begin{theorem} \label{shift}
There exists a unique, up to scaling, differential
  operator $S(m,\mu)$ (called the shift operator) 
of order $\sum_s m_s$, such that 
$S(m,\mu)H(\mu)=H(m+\mu)S(m,\mu)$. One has
  $S(m,\mu)L_q(\mu)=L_q(m+\mu)
S(m,\mu)$ for 
$W$-invariant polynomials $q$. 
The operator $S(m,\mu)$ is of degree $0$, 
has polynomial coefficients, and has symbol proportional to 
$\delta_m(x)\delta_m(\xi)$. \qed
\end{theorem}

\noindent
{\bf Example.} Let $W=\Z/2$ acting on $V=\R$ by $x\to -x$. 
In this case, there is only one number $m$, and one has:
$S(1,\mu)=x\partial-(2\mu+1)$,  
$$
S(m,\mu)=S(1,\mu+m-1)...S(1,\mu)=
c(x\partial-(2\mu+2m-1))...(x\partial-(2\mu+1)).
$$

The connection between the shift operator and the $\psi$-function 
is given by the following theorem, in which
`$:\;:$', the normal ordering
sign,  means that $x$ stands to the left of $\partial$.
 
\begin{theorem} \label{psishift} (see \cite{VSC}) 
One has $\psi(k,x)=S(m,0)^{(x)}e^{(k,x)}$. 
In other words, we have $S(m,0)=:P(x,\partial):$.\qed
\end{theorem}

The theorem follows from Theorem \ref{shift} and the fact that 
the function $\psi(k,x)$ of the form 
$P(k,x)e^{(k,x)}$ is uniquely determined already by the equations 
$L_q^{(x)}\psi=q(k)\psi$ for {\bf invariant} polynomials
$q$. 

\begin{corollary} $S(m,0)(\C[V])=Q_m$. 
Thus, $P_{Q_m}(t)+P_{\text{Ker}S(m,0)}(t)=(1-t)^{-n}$. \qed
\end{corollary} 

\footnotesize{

}
%\vskip 1cm

\footnotesize{
{\bf P.E.}: Department of Mathematics, Rm 2-165, MIT,
77 Mass. Ave, Cambridge, MA 02139\\ 
\hphantom{x}\quad\, {\tt etingof@math.mit.edu}

{\bf V.G.}: Department of Mathematics, University of Chicago, 
Chicago IL
60637, USA;\\ 
\hphantom{x}\quad\, {\tt ginzburg@math.uchicago.edu}}

\end{document}